\documentclass[twoside, 11pt]{article}

\usepackage{mathrsfs,amsfonts,amsmath}
\RequirePackage{ifthen,calc}
\usepackage{theorem}
\usepackage[dvips]{color}

\def\E{{\mathbb E}}
\def\P{{\mathbb P}}

 \catcode`@=11
 \def\@evenhead{\hbox to\textwidth{\footnotesize\rm\thepage \hfill
  {\it }}} 

 \def\@oddhead{\hbox to \textwidth{\footnotesize{\it
 Adaptive Elastic Net Method } \hfill\thepage}}

 \renewcommand{\section}{\makeatletter
 \renewcommand{\@seccntformat}[1]{{\csname the##1\endcsname.}\hspace{0.45em}}
 \makeatother \@startsection
{section}
{1}
{0pt}
{\baselineskip}
{0.5\baselineskip}
{\normalsize\bfseries\mathversion{bold}}}

\renewcommand{\subsection}{\makeatletter
 \renewcommand{\@seccntformat}[1]{{\csname the##1\endcsname.}\hspace{0.45em}}
 \makeatother \@startsection
{subsection}
{1}
{0pt}
{\baselineskip}
{0.5\baselineskip}
{\normalsize\bfseries\mathversion{bold}}}

\catcode`@=12

\newtheorem{theorem}{\noindent Theorem}[section]

\theoremheaderfont{\normalfont\bfseries}
\theorembodyfont{\slshape} {\theorembodyfont{\rmfamily}} {\theorembodyfont{\rmfamily}
\newtheorem{defn}{\noindent Definition}[section]}
{\theorembodyfont{\rmfamily} } {\theorembodyfont{\rmfamily}
}
{\theorembodyfont{\rmfamily} } {\theorembodyfont{\rmfamily}
}
{\theorembodyfont{\rmfamily} } {\theorembodyfont{\rmfamily} } {\theorembodyfont{\rmfamily}
}

 \def\beqlb{\begin{eqnarray}}\def\eeqlb{\end{eqnarray}}
 \def\beqnn{\begin{eqnarray*}}\def\eeqnn{\end{eqnarray*}}
 \numberwithin{equation}{section}
\def\2R{\mathbb{R}_+\times\mathbb{R}}
\def\qed{\hfill$\square$\smallskip}

\def\3R{\mathbb{R}_+\times\mathbb{R}_-}

\begin{document}
\title{\bf  Adaptive Elastic Net Method for Cox Model
\footnotetext{\hspace{-5ex}
${[1]}$ School of Mathematics, Guangxi University, Nanning 530004, China.
\\${[2]}$ School of Statistics, Shandong University of Finance and Economics,
Jinan 250014,  China.
\\${[3]}$ School of Mathematics and Statistics, Hechi University, Yizhou  546300, China
\\${[4]}$  E-mail:mathdsh@gmail.com
\newline
}}
\author{\small Chunhong Li$^1$, Xinxing Wei$^{1,3}$, Hongshuai Dai$^{2,4}$ }
\date{}

\maketitle

\begin{abstract}
In this paper, we study the Adaptive Elastic Net method for the Cox model. We prove the grouping effect and  oracle property of its estimators. Finally, we show these two properties by  an  empirical analysis and a numerical simulation, respectively.
\end{abstract}

\vspace{2mm}

\section{Introduction}
The aim of survival analysis is usually to identify  risk factors and their risk contributions. Often, many covariates are collected and  then a large parametric model is built. Hence, to efficiently select a subset of significant variables upon which the hazard function depends becomes an important and challenging task. Recently, many scholars used the variable selection techniques in studying linear regression models to deal with this kind of problems. For more information on this, see Fan and Li \cite{FL2002} and the references therein. However the treatment of using these methods directly also causes some problems. To overcome these drawbacks, statisticians have recently proposed a family of penalized partial likelihood methods to study the survival data, such as the Lasso method and so on.

The Lasso method introduced by Tibshirani \cite{T1997} is a penalized least squares method imposing a penalty on the regression coefficients. Due to the nature of the penalty, the Lasso method does both continuous shrinkage and automatic variable selection simultaneously. However, the Lasso estimator does not possess the oracle property and instability with high-dimensional data. Hence, Zou \cite{Z2006} proposed the Adaptive Lasso method, which has the oracle property. {\em Namely, the true regression coefficients that are zero are automatically estimated as zero, and the remaining coefficients are estimated as well as if the correct submodel were known in advance.} Contrast to the Lasso and Adaptive Lasso, the Elastic Net method proposed by Zou and Hastie \cite{ZH2005} is particularly useful when the number of predictors   is much bigger than that of observations. In addition, the Elastic Net method encourages a grouping effect, which means that {\em strongly correlated predictors tend to be in or out of the model together}. However, Fan and Li \cite{FL2001,FL2002} stated that the estimator of the Elastic Net method does not have the oracle property. Hence, Zou and Zhang \cite{ZZ2009} proposed the Adaptive Elastic Net method, which has the oracle property and grouping effect.

The Cox model\cite{C1972} is a classical method to deal with survival data. It is well-known that the Elastic Net method for the Cox model has grouping effect and can deal with the highly correlated data.  Inspired by these facts and Zou and Zhang \cite{ZZ2009}, in this paper, we study the Adaptive Elastic Net method for the Cox model and show that it has the grouping effect and oracle property.

The rest of this paper is organized as follows. In Section 2, we  introduce the Adaptive Elastic Net method for the Cox model. Section 3 is devoted to studying the grouping effect. The oracle property is discussed in Section 4. In Section 5, we show these two properties by an empirical analysis and a numerical simulation.

\section{Adaptive Elastic Net method}\label{sec2}
In this section, we  give the definition of the Adaptive Elastic Net method for the Cox model. We first recall some known facts about the Cox model. Recall that the hazard function for  an individual at the failure time $t$  is
\beqlb\label{a-23}
h(t)=h_0(t)\exp\{\beta^TX\},
\eeqlb
where $h_0(t)$ is a baseline hazard function,  $\beta=(\beta_1,\cdots,\beta_p)^T$ is the regression vector of  unknown coefficients, $X$ is the covariate of an individual. Let $R_i$ denote the risk set at time $t_i-0$, that is the set of  individuals who have not failed or been censored by that time. Furthermore, let $X_j=(x_{j1},\cdots,x_{jp})^T$ denote the value of $X$ for the $j$th individual and $X_i$  the value for the individual failing at time $t_i$. Suppose a random sample of $n$ individuals is chosen, then the likelihood function for inference about $\beta$ is given by:
\beqnn
L(\beta)=\Pi_{i=1}^n \frac{\exp\big(\sum_{k=1}^{p}\beta_kx_{ik}\big)}{\sum_{j\in R_i}\exp\big(\sum_{k=1}^p\beta_kx_{jk}\big)}.
\eeqnn
Therefore, the log-likelihood function is
\beqlb\label{1}
l(\beta)=\sum_{i=1}^n\Big\{\sum_{k=1}^p\beta_k x_{ik}-\ln\Big[\sum_{j\in R_i}\exp(\sum_{k=1}^p\beta_k x_{jk})\Big]\Big\}.
\eeqlb
By maximizing \eqref{1}, we can get the estimator of $\beta$.

From Tibshirani \cite{T1997}, and Fan and Li \cite{FL2002}, by minimizing the opposite number of \eqref{1} first, and then adding the appropriate penalty, we can get the Elastic Net estimator for the Cox model:
\beqlb\label{2}
\hat{\beta}_{(EN)}=\arg\min\bigg\{\frac{1}{n}\sum_{i=1}^n\Big\{-\beta^TX_i+\ln \Big[\sum_{j\in R_i}\exp\big(\beta^TX_j\big)\Big]\Big\}+\lambda_1\|\beta\|_1+\lambda_2\|\beta\|^2\bigg\},
\eeqlb
where $\lambda_1\geq 0$  and  $\lambda_2\geq 0$ are regularization parameters, and
\beqnn
\|\beta\|^2=\sum_{j=1}^p\beta_j^2,\;\textrm{and}\;\|\beta\|_1=\sum_{j=1}^p|\beta_j|.
\eeqnn
Moreover,
\eqref{2} can be rewritten as
\beqlb\label{3}
\hat{\beta}_{(EN)}=&&\arg\min\bigg\{\frac{1}{n}\sum_{i=1}^n\Big\{-\sum_{k=1}^p\beta_kx_{ik}+\ln \Big[\sum_{j\in R_i}\exp\big(\sum_{k=1}^p\beta_kx_{jk}\big)\Big]\Big\}\nonumber
\\&&\quad+\lambda_1\sum_{k=1}^p|\beta_k|+\lambda_2\sum_{k=1}^p\beta_k^2\bigg\}.
\eeqlb
Following Zou and Zhang \cite{ZZ2009}, we introduce the Adaptive Elastic Net estimator for the Cox model as follows.
\begin{defn}
The Adaptive Elastic Net estimator $\hat{\beta}_{(AEN)}=\Big(\hat{\beta}_{(AEN)_1},\cdots, \hat{\beta}_{(AEN)_p}\Big)^T$ for the Cox model  is defined  by
\beqlb\label{4}
\hat{\beta}_{(AEN)}=&&\arg\min\bigg\{\frac{1}{n}\sum_{i=1}^n\Big\{-\sum_{k=1}^p\beta_kx_{ik}+\ln \Big[\sum_{j\in R_i}\exp\big(\sum_{k=1}^p\beta_kx_{jk}\big)\Big]\Big\}\nonumber
\\&&\quad+\lambda^*_1\sum_{k=1}^p\hat{\omega}_k|\beta_k|+\lambda_2\sum_{k=1}^p\beta_k^2\bigg\},
\eeqlb
where $\hat{\omega}_k=(\big|\hat{\beta}_{(EN)_k}\big|)^{-\gamma}$ with $\gamma>0$.
\end{defn}
\section{Grouping effect}
In this section, we study the grouping effect of the Adaptive Elastic Net method for the Cox model. Before we state the main result, we need the following notation. Let  $x_a=(x_{1a},\cdots,x_{na})$  and $x_b=(x_{1b},\cdots,x_{nb})$, $a,b=1,\cdots,p$,  be highly correlated, and $\hat{\beta}_a(\lambda_1^*,\lambda_2)$  and $\hat{\beta}_b(\lambda_1^*,\lambda_2)$  denote  the estimators of  the $a$th variable $\tilde{X}_a$ and the $b$th variable $\tilde{X}_b$ in the covariate, respectively.  Following the notation in Andersen and Gill \cite{AG1982},  let $Y_i(t)=I(T_i\geq t, C_i\geq t)$, $N_i(t)= I(T_i\leq t, T_i\leq C_i)$, $\delta_i=I(T_i\leq C_i)$ and $Z_i=\min\{T_i,C_i\}$. Then,

\begin{theorem}\label{lem-2}
For the Cox model, given the data $(Z_i,\delta_i,X_i)$ and parameter $(\lambda_1^*,\;\lambda_2)$  , the responses are centered and the predictors are standardized. Let $\hat{\beta}(\lambda_1^*,\lambda_2)$  be the Adaptive Elastic Net estimator. Suppose that
\beqlb\label{a-1}
\hat{\beta}_a(\lambda_1^*,\lambda_2)\hat{\beta}_b(\lambda_1^*,\lambda_2)>0.
\eeqlb
Define
\beqnn
D_{\lambda_1^*,\lambda_2}(a,b)=\Big|\hat{\beta}_a(\lambda_1^*,\lambda_2)-\hat{\beta}_b(\lambda_1^*,\lambda_2)\Big|,
\eeqnn
then
\beqnn
D_{\lambda_1^*,\lambda_2}(a,b)\to 0,
\eeqnn
which means that $D_{\lambda_1^*,\lambda_2}(a,b)$ approximates to $0$.
\end{theorem}
{\it Proof:} By \eqref{a-1}, we have
\beqlb\label{a-2}
\textrm{sgn}\big\{\hat{\beta}_a(\lambda_1^*,\lambda_2)\big\}=\textrm{sgn}\big\{\hat{\beta}_b(\lambda_1^*,\lambda_2)\big\}
\eeqlb
and
\beqnn
\hat{\beta}_a(\lambda_1^*,\lambda_2)\neq 0\;\textrm{and}\;\hat{\beta}_b(\lambda_1^*,\lambda_2)\neq 0.
\eeqnn

Now, let $\hat{\beta}_m(\lambda_1^*,\lambda_2)\neq 0$ and at the point $\hat{\beta}(\lambda_1^*,\lambda_2)$,
$$
\frac{\partial L(\lambda_1^*,\lambda_2,\beta)}{\partial \beta_m}=0,
$$
where
\beqnn
L(\lambda_1^*,\lambda_2,\beta)=&&\frac{1}{n}\sum_{i=1}^n\Big\{-\sum_{k=1}^p\beta_kx_{ik}+\ln \Big[\sum_{j\in R_i}\exp\big(\sum_{k=1}^p\beta_kx_{jk}\big)\Big]\Big\}\nonumber
\\&&\quad+\lambda^*_1\sum_{k=1}^p\hat{\omega}_k|\beta_k|+\lambda_2\sum_{k=1}^p\beta_k^2.
\eeqnn
Then,
\beqnn
&&-\frac{1}{n}\sum_{i=1}^nx_{ia}+\frac{1}{n}\sum_{i=1}^n \frac{\sum_{j\in R_i}x_{ja}\exp\big(\sum_{k=1}^p\beta_kx_{jk}\big)}{\sum_{j\in R_i}\exp \big(\sum_{k=1}^p\beta_kx_{jk}\big)}
\\&&\qquad\qquad+\lambda_1^*\hat{\omega}_a\textrm{sgn}\big\{\hat{\beta}_a(\lambda_1^*,\lambda_2)\big\}+2\lambda_2\hat{\beta}_a(\lambda_1^*,\lambda_2)=0.
\eeqnn
Therefore,
\beqlb\label{5}
&&\hat{\beta}_a(\lambda_1^*,\lambda_2)=\nonumber
\\&&\quad\frac{1}{2\lambda_2}\bigg\{\frac{1}{n}\sum_{i=1}^nx_{ia}-\frac{1}{n}\sum_{i=1}^n \frac{\sum_{j\in R_i}x_{ja}\exp\big(\sum_{k=1}^p\beta_kx_{jk}\big)}{\sum_{j\in R_i}\exp \big(\sum_{k=1}^p\beta_kx_{jk}\big)}-\lambda_1^*\hat{\omega}_a\textrm{sgn}\big\{\hat{\beta}_a(\lambda_1^*,\lambda_2)\big\}\bigg\}.\;
\eeqlb
Similar to \eqref{5}, we get
\beqlb\label{6}
&&\hat{\beta}_b(\lambda_1^*,\lambda_2)=\nonumber
\\&&\quad\frac{1}{2\lambda_2}\bigg\{\frac{1}{n}\sum_{i=1}^nx_{ib}-\frac{1}{n}\sum_{i=1}^n \frac{\sum_{j\in R_i}x_{jb}\exp\big(\sum_{k=1}^p\beta_kx_{jk}\big)}{\sum_{j\in R_i}\exp \big(\sum_{k=1}^p\beta_kx_{jk}\big)}-\lambda_1^*\hat{\omega}_b\textrm{sgn}\big\{\hat{\beta}_b(\lambda_1^*,\lambda_2)\big\}\bigg\}.\;
\eeqlb
It follows from \eqref{a-2}, \eqref{5} and \eqref{6} that

\begin{align}\label{7}
&&\hat{\beta}_a(\lambda_1^*,\lambda_2)-\hat{\beta}_b(\lambda_1^*,\lambda_2)=\frac{1}{2\lambda_2}\Bigg\{\frac{1}{n}\sum_{i=1}^n\bigg[x_{ia}-x_{ib}+\frac{\sum_{j\in R_i}x_{jb}\exp\big(\sum_{k=1}^p\beta_kx_{jk}\big)}{\sum_{j\in R_i}\exp \big(\sum_{k=1}^p\beta_kx_{jk}\big)}\nonumber
\\&&-\frac{\sum_{j\in R_i}x_{ja}\exp\big(\sum_{k=1}^p\beta_kx_{jk}\big)}{\sum_{j\in R_i}\exp \big(\sum_{k=1}^p\beta_kx_{jk}\big)}
\bigg]+\lambda_1^*\textrm{sgn}\big\{\hat{\beta}_b(\lambda_1^*,\lambda_2)\big\}(\hat{\omega}_b-\hat{\omega}_a)\Bigg\}.
\end{align}
From Schoenfeld \cite{S1982}, we have
\begin{align}\label{8}
\hat{r}_{ir}:=x_{ir}-\E(x_{ir}|R_i)=x_{ir}-\frac{\sum_{l\in R_i}x_{lr}\exp(x_l\beta)}{\sum_{l\in R_i}\exp(x_l\beta)},
\end{align}
where $r=1,2,\cdots,p.$  Therefore, \eqref{7} is equivalent to
\beqnn
\hat{\beta}_a(\lambda_1^*,\lambda_2)-\hat{\beta}_b(\lambda_1^*,\lambda_2)=\frac{1}{2n\lambda_2}
\sum_{i=1}^n(\hat{r}_{ia}-\hat{r}_{ib})+\frac{\lambda_1^*}{2\lambda_2}\textrm{sgn}\big\{\hat{\beta}_b(\lambda_1^*,\lambda_2)\big\}
(\hat{\omega}_b-\hat{\omega}_a).
\eeqnn
Hence,
\beqlb\label{9}
\Big|\hat{\beta}_a(\lambda_1^*,\lambda_2)-\hat{\beta}_b(\lambda_1^*,\lambda_2)\Big|\leq \frac{1}{2n\lambda_2}
\sum_{i=1}^n\big|\hat{r}_{ia}-\hat{r}_{ib}\big|+\frac{\lambda_1^*}{2\lambda_2}
\big|\hat{\omega}_b-\hat{\omega}_a\big|.
\eeqlb
Since $x_a$  and $x_b$  are highly correlated, i.e.,
\beqnn
\E\big[x_ax_b^T\big]\to 1.
\eeqnn
 Then, we have for an individual $i$
\beqnn
|x_{ia}-x_{ib}|\to 0,\;\textrm{and}\;\big|\E(x_{ia})-\E(x_{ib})\big|\to 0.
\eeqnn

Hence,
\beqlb\label{a-3}
\Big|\big[x_{ia}-\E(x_{ia}|R_i)\big]-\big[x_{ib}-\E(x_{ib}|R_i)\big]\Big|\to 0.
\eeqlb
By \eqref{8} and \eqref{a-3}, we have
\beqlb\label{a-5}
|\hat{r}_{ia}-\hat{r}_{ib}|\to 0\;\textrm{and}\;\big|\E(\hat{r}_{ia})-\E(\hat{r}_{ib})\big|\to 0.
\eeqlb

Since
\beqnn
\hat{\omega}_a=\Big(\big|\hat{\beta}_{(EN)a}\big|\Big)^{-\gamma}\;\textrm{and}\;\hat{\omega}_b=\Big(\big|\hat{\beta}_{(EN)b}\big|\Big)^{-\gamma},
\eeqnn
we have
\beqlb\label{a-4}
\big|\hat{\omega}_a-\hat{\omega}_b\big|\to 0.
\eeqlb
It follows from \eqref{9}, \eqref{a-5} and \eqref{a-4} that
\beqnn
D_{\lambda_1^*,\lambda_2}(a,b)\to 0.
\eeqnn
The lemma holds.
\qed

\section{Oracle property }
In this section, we study the oracle property of the Adaptive Elastic Net method for the time-dependent Cox model. To emphasize the time dependence, we rewrite the model \eqref{a-23} as follows
\beqlb\label{a-6}
h(t|X)=h_0(t)\exp\big\{\beta^TX(t)\big\},
\eeqlb
where the covariate  $X(t)$ is time-dependent. Hence, the Adaptive Elastic Net estimator for this Cox model is:
\beqlb\label{a-7}
\hat{\beta}_{(AEN)}=\arg \min\bigg\{-\frac{1}{n}l_n(\beta)+\lambda^*_1\sum_{k=1}^p\hat{\omega}_k|\beta_k|
+\lambda_2\sum_{k=1}^p\beta_k^2\bigg\},
\eeqlb
where $l_n(\beta)$ is the partial log-likelihood function.

Next we consider
 \beqnn
 Q_n(\beta)=l_n(\beta)-n\lambda^*_1\sum_{k=1}^p\hat{\omega}_k|\beta_k|
-n\lambda_2\sum_{k=1}^p\beta_k^2.
 \eeqnn

We suppose that the real  parameter $\beta_0=(\beta_{01},\cdots,\beta_{0p})^T$ is sparse, and $A=\big\{k:\beta_{0k}\neq 0\big\}=\{1,\cdots,p_0\}$ with $p_0<p$. Hence $\beta_0=\big(\beta_{A},\beta_{A^c}\big)^T$  with $\beta_{A^c}=\bf{0}$, where $\bf{0}$ is the zero vector. Hence the corresponding estimator $\hat{\beta}_0$  takes the form of  \beqnn \hat{\beta}_0=\big(\hat{\beta}_{A},\hat{\beta}_{A^c}\big)^T . \eeqnn
 Let $I(\beta_0)$ be the Fisher information matrix. We also introduce the following notation. For any matrix $B$,
 \beqnn
 \|B\|=\sup_{ij}|b_{ij}|,
 \eeqnn
 and for any vector $a$
 \beqnn
 a^{\otimes 0}=1,\;a^{\otimes 1}=a,\;a^{\otimes 2}=a\otimes a,\;\|a\|=\sup_i|a_i|,\;\textrm{and}\;|a|=\big(\sum a_i^2\big)^\frac{1}{2},
 \eeqnn
where $a\otimes b$ is the $p\times p $ matrix $ab^T$  for any vectors $a,b \in\mathbb{R}^p$.
 Before we state our main result in this section, similar to Fan and Li \cite{FL2002}, we need the following conditions:
 \begin{itemize}
 \item[(a)]
 \beqlb\label{a-8}
 \int_0^1h_0(t)dt<\infty, n\lambda_1^*\to \infty, \textrm{and}\;\frac{\lambda_2}{\sqrt{n}}\to 0,\;\textrm{as}\; n\to \infty;
 \eeqlb
 \item[(b)]There exists a neighborhood $\Omega$ of $\beta_0$ such that for any $\beta\in\Omega$,
 \beqlb\label{a-9}
 \E\big[\sup_{t\in[0,\;1],\beta\in\Omega}Y(t)X^T(t)X(t)\exp\big\{\beta^TX(t)\big\}\big]<\infty;
 \eeqlb
 Moreover,
 \item[(c)] for any $\beta\in\Omega$, we have
 \beqnn
 s_0(\beta,t)&&=\E\big[Y(t)\exp\big\{\beta^TX(t)\big\}\big],\\
  s_1(\beta,t)&&=\E\big[Y(t)X(t)\exp\big\{\beta^TX(t)\big\}\big],\;\textrm{and}\;
  \\ s_2(\beta,t)&&=\E\big[Y(t)X(t)X^T(t)\exp\big\{\beta^TX(t)\big\}\big],
 \eeqnn
 where $s_i(\beta,t),i=0,1,2,$ satisfies the following:
 \begin{itemize}
 \item[(I)]$s_0(\beta,t), s_1(\beta,t)$ and $s_2(\beta,t)$ are uniformly continuous in $t\in[0,\;1]$ for $\beta\in\Omega$.
 \item[(II)]$s_2(\beta, t)$  and $s_1(\beta,t)$ are bounded on $\Omega\times[0,\;1]$. Moreover assume $s_0(\beta,t)$ is positive and bounded away from zero on $\Omega\times[0,\;1]$.
 \end{itemize}
 \item[(d)] Define:
     \beqnn
 v(\beta,t)=\frac{s_2(\beta,t)}{s_0(\beta,t)}-e^{\otimes2}
     \eeqnn
 with $$
 e=\frac{s_1(\beta,t)}{s_0(\beta,t)}.
 $$
 Then,
 the fisher information matrix
  \beqnn
  I(\beta_0)=\int_0^1v(\beta_0,t)s_0(\beta_0,t)h_0(t)dt
  \eeqnn
 \end{itemize}
 is finite positive definite.

In addition, the regularity conditions  in Anderson and Gill \cite{AG1982} are assumed in the whole section. Then, we have the following theorem.

\begin{theorem}\label{thm-1}
If  (a)-(d) hold, then the Adaptive Elastic Net estimator  $\hat{\beta}$ for the time-dependent Cox model has the sparsity, i.e.,
\beqnn
\P\big(\hat{\beta}_{A^c}={\bf 0}\big)\to 1.
\eeqnn
\end{theorem}

{\it Proof:}
We note that the partial log-likelihood function of the time-dependent Cox model is:
\beqlb\label{11}
l_n(\beta)=\sum_{i=1}^n\int_0^1\beta^TX_i(t)dN_i(t)-\int_0^1\log\big[\sum_{i=1}^nY_i(t)\exp\big\{\beta^TX_i(t)\}\big]d\tilde{N}(t)
\eeqlb
where  $\tilde{N}(t)=\sum_{i=1}^nN_i(t)$.

Then, for each $\beta$ in the neighborhood $\Omega$ of  $\beta_0$, we have
\beqlb\label{12}
\frac{1}{n}\{l_n(\beta)-l_n(\beta_0)\}=f(\beta)+O_p\Big(\frac{\|\beta-\beta_0\|}{\sqrt{n}}\Big),
\eeqlb
where $O_p(\cdot)$ denotes convergence rate, and
\beqlb\label{a-11}
f(\beta)=\int_0^1\Big[(\beta-\beta_0)^Ts_1(\beta_0,t)-\log \big\{\frac{s_0(\beta,t)}{s_0(\beta_0,t)}\}s_0(\beta_0,t)\Big]h_0(t)dt.
\eeqlb

Let $\beta=\beta_0+\frac{u}{\sqrt{n}}$, where $\|u\|\leq C$ for some large enough constant $C$. According to the theorem 3.1 in Fan and Li \cite{FL2001}, we know that for any $\epsilon>0$,
\beqnn
\P\Big\{\sup_{\|u\|=C}Q_n(\beta)<Q_n(\beta_0)\Big\}\geq 1-\epsilon.
\eeqnn
Then
\beqlb\label{13}
\Big\|\hat{\beta}_0-\beta_0\Big\|=O_p(n^{-\frac{1}{2}}).
\eeqlb
By the Taylor expansion of $l_n(\beta)$, we have
\beqlb\label{a-10}
l_n(\beta)&&=l_n(\beta_0)+l'_n(\beta_0)(\beta-\beta_0)+O_p(\sqrt{n}\|\beta-\beta_0\|).
\eeqlb
By \eqref{12} and \eqref{a-10}, we obtain
\beqlb\label{14}
l_n(\beta)&&=l_n(\beta_0)+nf(\beta)+O_p(\sqrt{n}\|\beta-\beta_0\|),
\eeqlb
where
$
f(\beta)
$
is given by \eqref{a-11}.

On the other hand, by the condition (c), we have
\beqnn
\frac{\partial f(\beta)}{\partial \beta}=\int_0^1\Big[\frac{s_1(\beta_0,t)}{s_0(\beta_0,t)}-\frac{s_1(\beta,t)}{s_0(\beta,t)}\Big]s_0(\beta_0,t)h_0(t)dt,
\eeqnn
and
\beqnn
-\frac{\partial^2 f(\beta)}{\partial \beta\partial \beta^T}=\int_0^1\Big[\frac{s_2(\beta,t)s_0(\beta,t)-s_1(\beta,t)s^T_1(\beta,t)}{\big[s_0(\beta,t)\big]^2}\Big]s_0(\beta_0,t)h_0(t)dt.
\eeqnn
Hence
\beqlb\label{a-12}
f(\beta_0)=0,\;\frac{\partial f(\beta)}{\partial \beta}|_{\beta=\beta_0}={\bf 0},\;\textrm{and}\;-\frac{\partial^2 f(\beta)}{\partial \beta\partial \beta^T}=I(\beta_0).
\eeqlb
From \eqref{a-12} and the Taylor expansion, we have
\beqlb\label{a-13}
f(\beta)=-\frac{1}{2}(\beta-\beta_0)^T\Big\{I(\beta_0)+o(1)\Big\}(\beta-\beta_0).
\eeqlb
Next we study \eqref{14}. Since
\beqnn
\frac{\partial l_n(\beta)}{\partial \beta_k}=O_p(\sqrt{n}),\;k=p_0+1,\cdots,p,
\eeqnn
we have
\beqlb\label{15}
\frac{\partial Q_n(\beta)}{\partial \beta_k}=O_p(\sqrt{n})-n\lambda_1^*\sum_{k}\omega_k\textrm{sgn}\{\beta_k\}-2n\lambda_2\sum_{k}\beta_k.
\eeqlb
Noting
\beqnn
n^{\frac{3}{2}}\big(|\beta_k|^{\gamma}-0\big)=O_p(1),
\eeqnn
we obtain
\beqlb\label{a-14}
\frac{\partial Q_n(\beta)}{\partial \beta_k}=n^{\frac{3}{2}}\Big\{O_p(n^{-1})-n\lambda_1^*\frac{\textrm{sgn}\{\beta_k\}}{O_p(1)}-2\frac{\lambda_2}{\sqrt{n}}\sum_{k}\beta_k\Big\}.
\eeqlb
Combing \eqref{a-8} and \eqref{a-14}, we get that the sign of $\frac{\partial Q_n(\beta)}{\partial \beta_k}$ is determined only by $\beta_k$, where $$\beta_k\in\big(-Cn^{-\frac{1}{2}},\;Cn^{-\frac{1}{2}}\big),\;k=p_0+1,\cdots,p.$$
Therefore,
\beqlb\label{16}
\big\|\hat{\beta}_A-\beta_A\big\|=O_p(n^{-\frac{1}{2}}).
\eeqlb
By \eqref{13} and \eqref{16},
\beqnn
Q_n(\beta_A,0)=\max_{\|\beta_{A^c}\|\leq C n^{-\frac{1}{2}}}Q_n(\hat{\beta}_A,\hat{\beta}_{A^c}).
\eeqnn
Therefore
\beqnn
\P\big(\hat{\beta}_{A^c}={\bf 0}\big)\to 1.
\eeqnn
The proof is completed.
\qed

Next we study its asymptotic normality.
\begin{theorem}
Suppose that the conditions (a)-(d) hold. Then the Adaptive Elastic Net estimator $\hat{\beta}$ for the time-dependent Cox model has the following asymptotic normality:
\beqnn
\sqrt{n}\big(\hat{\beta}_A-\beta_A)\stackrel{D}{\to} N({\bf 0},\;I_1^{-1}(\beta_A)),
\eeqnn
where $\stackrel{D}{\to}$ denotes convergence in distributions.
\end{theorem}
\noindent{\bf Proof:} According to the proof of Theorem \ref{thm-1}, we know that there exists $\hat{\beta}_A$ such that for $k=1,\cdots,p_0$
\begin{align}\label{a-15}
\frac{\partial Q_n(\beta)}{\partial \beta_k}\Big|_{\beta=(\hat{\beta}_A,{\bf 0})^T}=0.
\end{align}
Let $U_n(\beta)$ be the score function of $l_n(\beta)$, that is
\begin{align}\label{a-16}
U_n(\beta)=\sum_{i=1}^n\int_0^1X_i(t)dN_i(t)
-\int_0^1\frac{\sum_{i=1}^nY_i(t)X_i(t)\exp\{\beta^TX_i(t)\}}
{\sum_{i=1}^nY_i(t)\exp\{\beta^TX_i(t)\}}d\tilde{N}(t).
\end{align}
Moreover, define
\beqlb\label{a-17}
\hat{I}(\beta)&&=\int_0^1\Bigg[\frac{\sum_{i=1}^nY_i(t)X_i(t)X_i^T(t)\exp\{\beta^TX_i(t)\}}{\sum_{i=1}^nY_i(t)\exp\{\beta^TX_i(t)\}}\nonumber
\\&&-\frac{[\sum_{i=1}^{n}Y_i(t)X_i(t)\exp\{\beta^TX_i(t)\}][\sum_{i=1}^n Y_i(t)X_i(t)\exp\big\{\beta^TX_i(t)\}]^T}{\big[\sum_{i=1}^nY_i(t)\exp\big\{\beta^TX_i(t)\big\}\big]^2}\Bigg]d\tilde{N}(t).\qquad
\eeqlb
On the other hand, we have
\beqlb\label{17}
\frac{\partial Q_n(\beta)}{\partial \beta_k}\Big|_{\beta=(\hat{\beta}_A,{\bf 0})^T}=\frac{\partial l_n(\beta)}{\partial \beta_k}\Big|_{\beta=(\hat{\beta}_A,{\bf 0})^T}-n\lambda_1^*\sum_k \omega_k\textrm{sgn}\{\beta_k\}-2n\lambda_2\sum_k\beta_k=0.
\eeqlb

According to the Taylor expansion, \eqref{17} is equivalent to
\beqlb\label{18}
U_A(\beta_0)-\hat{I}_1(\tilde{\beta})(\hat{\beta}_A-\beta_A)-n\lambda_1^*\sum_k\omega_k\textrm{sgn}\{\beta_k\}-2n\lambda_2\sum_k\beta_k=0,
\eeqlb
where $\tilde{\beta}=(\hat{\beta}_0,\beta_0)$.

According to \eqref{a-8} and  Andersen and Gill \cite{AG1982}, we have
\beqlb\label{a-18}
\frac{1}{\sqrt{n}}U_A(\beta_0)\to N\big({\bf 0},I_1(\beta_A)\big),
\eeqlb
and
\beqlb\label{19}
\frac{1}{n}\hat{I}_1(\tilde{\beta})\to I_1(\beta_A),
\eeqlb
where $U_A(\beta_0)$ is composed by the first $p_0$ elements of $U(\beta_0)$, and $I_1(\beta_A)$ is a $p_0\times p_0$ submatrix of $I(\beta_0).$

On the other hand, if we assume
\beqlb\label{a-21}
\sqrt{n}\lambda_1^*\to\lambda_0,
\eeqlb
then
\beqlb\label{20}
\frac{1}{\sqrt{n}}(\hat{\beta}_A-\beta_A)=\hat{I}_1^{-1}(\tilde{\beta})\Big[\frac{1}{\sqrt{n}}
U_A(\beta_0)-\sqrt{n}\lambda_1^*\sum_{k}\omega_k\textrm{sgn}\{\beta_k\}\Big]
+O_p(1).
\eeqlb
Note that from \eqref{19},
\beqlb\label{a-19}
n\hat{I}_1^{-1}(\tilde{\beta})\to
I_1^{-1}(\beta_A).
\eeqlb
Then \eqref{20} can be rewritten as:
\beqlb\label{21}
\sqrt{n}(\hat{\beta}_A-\beta_A)=I_1^{-1}(\beta_A)\Big[\frac{1}{\sqrt{n}}
U_A(\beta_0)-\lambda_0b_1\Big]+O_p(1),
\eeqlb
where $$b_1=\sum_k\omega_k\textrm{sgn}\{\beta_k\},\;k=1,\cdots,p_0.$$
According to the Slutsky's Theorem, we have as $n\to\infty$,
\beqlb\label{22}
\sqrt{n}(\hat{\beta}_A-\beta_A)\to N\Big(-\lambda_0I_1^{-1}(\beta_A)b_1,I_1^{-1}(\beta_A)\Big).
\eeqlb
Hence, if $\lambda_0=0$, then \eqref{21} and \eqref{22} can be simplified to
\beqlb\label{a-22}
\sqrt{n}(\hat{\beta}_A-\beta_A)=I_1^{-1}(\beta_A)\big[\frac{1}{\sqrt{n}}U_A(\beta_0)\big]+O_p(1),
\eeqlb
and
\beqlb\label{a-20}
\sqrt{n}(\hat{\beta}_A-\beta_A)\to N\Big({\bf 0},I_1^{-1}(\beta_A)\Big),\;\textrm{as}\; n\to\infty,
\eeqlb
respectively. The proof of the theorem is finished.\qed
\section{Empirical analysis}
Theorem 3.1 reveals that the Adaptive Elastic Net estimator for the Cox model enjoys the grouping effect. Next, we do an empirical analysis.

The data come from a questionnaire, which studies the mobile phone cards' usage of college students. $x_1,\cdots,x_9$ and $x_{10}$ denote Sex, Grade, Position, Nation, Registered residence, School address, Operators, The average monthly telephone charges, The quality of service and The average monthly living expenses, respectively. In addition, the study was recorded by years, and from 2007 to 2014. We finally got the 380 effective questionnaires. The training set number is 300, and the test set number is 80. Part of the data is listed in the following table 1.
\begin{table} [h]
\centering
\caption{ Data sample}
\label{T1}
\begin{tabular}{c|c|c|cccccccccc}
\hline
No. & Study times& Status&$x_1$&$x_2$&$x_3$&$x_4$&$x_5$&$x_6$&$x_7$&$x_8$&$x_9$&$x_{10}$ \\
1 & 6& censored &2&6&1&1&2&1&1&3&2&2\\
2 & 1& loss &2&3&1&1&1&1&1&2&2&1\\
3 & 1& loss &1&3&2&2&1&1&2&2&2&1\\
4 & 1& loss &1&6&1&1&1&1&2&2&1&1\\
5 & 4& censored &2&4&2&2&1&1&1&5&2&4\\
$\cdots$ & $\cdots$& $\cdots$&$\cdots$&$\cdots$&$\cdots$&$\cdots$&$\cdots$&$\cdots$&$\cdots$&$\cdots$&$\cdots$&$\cdots$
\\
380 & 2& loss &2&2&2&2&2&1&1&2&3&1\\
\hline
\end{tabular}
\end{table}

 Since $x_8$ and $x_{10}$ are medium correlation, we did the variable selection by the Lasso method, the Adaptive Lasso method(ALasso), the Elastic Net method(EN) and the Adaptive Elastic Net method(AEN), respectively. The selected variables are in Table \ref{T3}, and the coefficient estimators are  in Table \ref{T2}.
\begin{table}[h]
\centering
\caption{Statistics results for the real data}
\label{T3}
\begin{tabular}{ccccccccccc}
\hline
$\;$&$\;$&$\;$&$\;$&$\;$& Variables&  selected&  in the&  model&
\\
\hline
Lasso & $\;$&$x_2$&$\;$&$x_4$&$x_5$&$x_6$&$\;$&$\;$&$x_9$&$x_{10}$ \\
ALasso&  $\;$&$x_2$&$\;$&$\;$&$x_5$&$\;$&$\;$&$\;$&$x_9$&$x_{10}$\\
EN & $\;$&$x_2$&$\;$&$x_4$&$x_5$&$x_6$&$\;$&$x_8$&$x_9$&$x_{10}$\\
AEN & $\;$&$x_2$&$\;$&$\;$&$x_5$&$x_6$&$\;$&$x_8$&$x_9$&$x_{10}$\\
\hline
\end{tabular}
\end{table}

\begin{table}[h]
\centering
\caption{Coefficient estimators obtained by the real data}
\label{T2}
\begin{tabular}{ccccccccccc}
\hline
$\;$&$x_1$&$x_2$&$x_3$&$x_4$&$x_5$&$x_6$&$x_7$&$x_8$&$x_9$&$x_{10}$
\\
\hline
Lasso & $0$ &$-0.092$& $0.005$& $-0.011$ &  $0.125$&  $-0.246$ &$0$ &$0$&$0.277$    &$0.012$ \\
ALasso&  $0$&$-0.042$& $0$    & $0$      &  $0.096$&  $0$      &$0$ &$0$&$0.018$    &$0.047$\\
EN    & $0$ &$-0.112$& $0.004$& $-0.023$&   $0.124$& $-0.284$ &$0$ &$0.019$&$0.313$&$0.026$\\
AEN   & $0$ &$-0.120$& $0.002$& $0$      &   $0.125$&  $-0.301$&$0$ &$0.030$&$0.329$&$0.036$\\
\hline
\end{tabular}
\end{table}
From  Table \ref{T2}, we obtain the following:
\begin{itemize}
\item[(1)] These four methods do not select $x_1$ and $x_7$ into the model,  which shows that the respondents' gender and operators have no effect to the usages of mobile phone card.
\item[(2)]  $x_2$ is negative, which shows the senior students have a low  probability of loss.  $x_9$ is positive, which shows that the higher the students' average monthly living expenses is, the greater the loss probability is. Similarly, $x_{10}$ is positive, which means that  the worse the operator's customer service quality is, the greater the loss probability is. It is consistent with the actual situation.
\item[(3)] The coefficient estimators obtained by the AEN  are the most close to the true model.
\item[(4)] For $x_8$ and $x_{10}$, both the Lasso and ALasso select $x_8$ only, while the EN  and AEN  can select both into the model, suggesting that these two methods can select the all strongly correlated variables into the model, and their estimators of coefficients are almost the same, which reflects the grouping effect.
\end{itemize}

Theorems 3.1, 4.1 and 4.2 reveal that the Adaptive Elastic Net estimator for the Cox model enjoys the grouping effect and oracle property. Next, we show these properties through a numerical simulation.

Let $x_i\sim N(0,1), i=1,2,5,6,8,9,10$. Moreover, let $x_2=x_3,x_6=x_7$ and $x_4=2x_1+\frac{1}{2}x_2+\frac{1}{2}x_3$. Then $x_2$ and $x_3$  are strongly correlated, so as $x_6$ and $x_7$. Moreover, there exists the linear relationship between $x_1,x_2,x_3$ and $x_4$. Consider the following Cox model,
\beqnn
h(t)=h_0(t)\exp\big(\sum_{i=1}^{10} \beta_ix_i\big),
\eeqnn
where $t\sim U[0,\;1]$. The real parameter $\beta$  is $(-1,2,2,0,\frac{1}{2},1,1,0,0,0)^T$. Then we did the simulation with $n=1000$ and $p=10.$

We  used the Lasso, ALasso, EN and AEN to do variable selection, respectively. Let $\lambda_2=\frac{1}{3}$, $\gamma=3$, and the other parameters be selected by the cross validation method \cite{V1993}. By using the Lars algorithm \cite{EHJ2004}, we obtain the coefficient estimators. See Table \ref{T4}.

\begin{table}[h]
\centering
\caption{Coefficient estimators obtained by numerical simulation}
\label{T4}
\begin{tabular}{ccccccccccc}
\hline
$\;$&$x_1$&$x_2$&$x_3$&$x_4$&$x_5$
\\
\hline
Lasso & $-0.99553$ &$3.89255$& $0.09627$& $0$ &  $0.49304$  \\
ALasso&  $-0.99796$&$3.89967$& $0.09901$& $0$ &  $0.49936$ \\
EN    & $-0.99509$ &$1.99726$& $1.99726$& $0$ &   $0.50021$ \\
AEN   & $-0.99947$ &$1.99951$& $1.99951$& $0$ &   $0.49992$ \\
\hline
\\$\;$&$x_6$&$x_7$&$x_8$&$x_9$&$x_{10}$&\\
\hline
Lasso & $1.88214$& $0.05014$ &$0.00794$&$0.00261$ &$0.00385$ \\
ALasso& $1.90981$& $0.03938$ &$0$&$0$  &$0.00013$\\
EN    & $0.99854$& $0.99854$ &$0.00261$&$0.00016$&$0.00248$\\
AEN   & $0.99976$& $0.99976$ &$0$&$0.00010$&$0$\\
\hline
\end{tabular}
\end{table}

From  Table \ref{T4}, we get:
\begin{itemize}
\item[(1)] The coefficient estimators obtained by the AEN are the most close to the true model.
\item[(2)] None of the four methods selects $x_4$ into the model, which implies that these four methods are all able to deal with the collinearity problems.
\item[(3)] We look at grouped variables $(x_2, x_3)$ and  $(x_6,x_7)$. The AEN  and  EN  can select  all the strongly correlated variables into the model, and the coefficient estimators of these two groups are the same. While the Lasso  and  ALasso  select  $x_2$ and $x_6$, respectively. This indicates the AEN and  EN  enjoy the grouping effect.
\item[(4)] We focus on the variables $x_8$, $x_9$ and $x_{10}$. The ALasso  and  AEN  can get more accurate  estimators for the estimation of zero variables than the other two methods do.  This indicates that the AEN  has the oracle property.
\end{itemize}

\section{Conclusion}
In this paper, we study the Adaptive Elastic Net method for the Cox model. We show that it has the grouping effect and  oracle property. These two properties are showed by an empirical analysis and a numerical simulation. In these examples, the Adaptive Elastic Net  and  Elastic Net  can make up for the lack of the Lasso  and  Adaptive Lasso, and can select all the strongly correlated variables into the model, i.e., the Adaptive Elastic Net method for the Cox model enjoys the grouping effect. In addition, the Adaptive Elastic Net method for the Cox model has the oracle property.

\bigskip
\noindent{\bf Acknowledgments}\;    This work is supported by the  National Natural Science Foundation of China (No.11361007), the Guangxi Natural Science Foundation (Nos.2014GXNSFCA118001 and 2012GXNSFBA053010) and the Project for Fostering Distinguished Youth Scholars of Shandong University of Finance and Economics.


\end{document}